\documentclass{amsart}
\usepackage{amsmath}
\usepackage{eufrak}

\newtheorem{theorem}{\textbf{Theorem}}

\def\a {\alpha}

\def\Z {\mathbb{Z}}
\def\Q {\mathbb{Q}}
\def\R {\mathbb{R}}

\def\QQ {\overline{\Q}}
\def\C {\mathbb{C}}

\def\den{\mathop{\rm den}}

\theoremstyle{remark}

\numberwithin{equation}{section}

\begin{document}

\title[ON A QUESTION PROPOSED BY K. MAHLER]{ON A QUESTION PROPOSED BY K. MAHLER CONCERNING LIOUVILLE NUMBERS}
%    Information for first author
%    \thanks will become a 1st page footnote.

\author{JEAN LELIS}
\address{DEPARTAMENTO DE MATEM\'{A}TICA, UNIVERSIDADE DE BRAS\'ILIA, BRAS\'ILIA, DF, BRAZIL}
\email{jeancarlos@mat.unb.br}

\author{DIEGO MARQUES}
\address{DEPARTAMENTO DE MATEM\'{A}TICA, UNIVERSIDADE DE BRAS\'ILIA, BRAS\'ILIA, DF, BRAZIL}
\email{diego@mat.unb.br}

\author{JOSIMAR RAMIREZ}
\address{DEPARTAMENTO DE MATEM\'{A}TICA, UNIVERSIDADE DE BRAS\'ILIA, BRAS\'ILIA, DF, BRAZIL}
\email{josimar@mat.unb.br}

%    General info


\begin{abstract}
In 1906, Maillet proved that given a non-constant rational function $f$, with rational coefficients, if $\xi$ is a Liouville number, then so is $f(\xi)$. Motivated by this fact, in 1984,  Mahler raised the question about the existence of transcendental entire functions with this property. In this work, we provide an uncountable subset of Liouville numbers for which there exists a transcendental entire function taking this set into the set of the Liouville numbers.
\end{abstract}

\maketitle

%%  SECTION 1
\section{Introduction}

A \textit{transcendental function} is a function $f(x)$ such that the only complex polynomial satisfying $P(x, f(x)) =0$, for all $x$ in its domain, is the null polynomial. For instance, the trigonometric functions, the exponential function, and their inverses.

The study of the arithmetic behavior of transcendental functions at complex points has attracted the attention of many mathematicians for decades. The first result concerning this subject goes back to 1884, when Lindemann proved that the transcendental function $e^z$ assumes transcendental values at all nonzero algebraic point. In 1886, Weierstrass gave an example of
a transcendental entire function which takes rational values at all rational points. Later, St\"{a}ckel \cite{19} proved
that for each countable subset $\Sigma\subseteq \C$ and each dense subset $T\subseteq \C$, there is a transcendental entire
function $f$ such that $f(\Sigma) \subseteq T$ (F. Gramain showed that St\"{a}ckel's theorem is valid if $\Sigma$ and $T$ are subsets of $\R$). Another construction due to St\"{a}ckel \cite{20} produces a transcendental entire function $f$ whose
derivatives $f^{(t)}$, for $t = 0, 1, 2, \ldots$, all map $\QQ$ into $\QQ$. Two years later, G. Faber refined this result by showing the existence of a transcendental entire function such that $f ^{(t)}(\QQ)\subseteq  \Q(i)$, for all $t \geq 0$. In 1968, van der Poorten \cite{alf} proved the existence of a transcendental function $f$, such that $f^{(s)}(\alpha)\in \Q(\alpha)$, for all $\alpha\in \QQ$.

Note that all the previously mentioned results deal with the arithmetic behavior of a countable set by a transcendental function. 

A real number $\xi$ is called a \textit{Liouville number}, if there exist infinitely many rational numbers $(p_n/q_n)_{n\geq 1}$, with $q_n\geq 1$ and such that
\[
	0<\left| \xi - \frac{p_n}{q_n} \right| < \frac{1}{q_n^{n}}.
\]

It is well-known that the set of the Liouville numbers $\mathbb{L}$ is a $G_{\delta}$-dense set and therefore an uncountable set.

In his pioneering book, Maillet \cite[Chapitre III]{mai} discusses some arithmetic properties of Liouville numbers. One of them is that, given a non-constant rational function $f$, with rational coefficients, if $\xi$ is a Liouville number, then so is $f(\xi)$. Motivated by this fact, in 1984, as the first problem in his paper {\it Some suggestions for further research}, Mahler \cite{mahler1} raised the following question (this question also appeared in other texts as for example in the Bugeaud's book \cite[p. 215]{bugeaud} and in Waldschmidt's paper \cite[p. 281]{W})\\

\noindent
{\bf Mahler's question.} \textit{Are there transcendental entire functions $f(z)$ such that if $\xi$ is any Liouville number, then so is $f(\xi)$?}\\

He also said that: ``The difficulty of this problem lies of course in the fact that the set of all Liouville
numbers is non-enumerable". Bernik and Dombrovski\u i \cite{bernik} and Alnia\c cik \cite{AL} obtained some results related to this question. Also, recently, some authors (see \cite{gugu, josimar, jo}) constructed classes of Liouville numbers which are mapped into Liouville numbers by transcendental entire functions. 

We remark about the existence of more specific classes of Liouville numbers in the literature. For example, the strong and semi-strong Liouville numbers (see, for instance, \cite{Pet}). Here, we shall define an uncountable subclass of the strong Liouville numbers which we named as {\it ultra-strong Liouville numbers}: a real number $\xi$ is called an ultra-strong Liouville number, if the sequence of the convergents of its continued fraction satisfies
\[
0<\left| \xi - \frac{p_n}{q_n} \right| < \frac{1}{q_n^{n}},\ \mbox{for\ all}\ n\geq 1.
\]
We denote this set by $\mathfrak{L}$. Define a sequence $A=(a_n)_n$ by $a_1=a_2=a_3=1$ and $a_j\in \{v_{j-1},v_{j-1}+1\}$, for $j\geq 4$, where $v_{j-1}:=(\prod_{k=1}^{j-1}(a_k+1))^{j-3}$. Then the number $\xi_A:=[0,a_1,a_2,a_3,a_4,\ldots]$ is an ultra-strong Liouville number. In fact, if $[0,a_1,\ldots, a_n]=p_n/q_n$, then, by construction, $a_{n+1}>(\prod_{k=1}^{n}(a_k+1))^{n-2}>q_n^{n-2}$ (here, we used the well-known inequality $(a_1+1)\cdots (a_n+1)>q_n$). Thus
\[
0<\left|\xi_A-\frac{p_n}{q_n}\right|<\frac{1}{a_{n+1}q_n^2}<\frac{1}{q_n^n}
\]
as desired. The set $\mathfrak{L}$ is uncountable because there exists a binary tree of possibilities for $\xi_A$, since we have two possibilities for $a_k$ in each step ($k\geq 4)$.

In this paper, we prove the following result:

\begin{theorem}\label{teo1}
Let $(s_n)_{n\geq 1}$ be sequence of positive integers satisfying that, for any given $k\geq 1$, the quotient $s_n/s_{n-1}^k$ tends to infinity as $n\to \infty$. Let $F:\C\to \C$ be a function defined by 
\[
F(z)=\displaystyle\sum_{k\geq 1}\frac{\a_k}{10^{k!}}z^k,
\]
where $\a_k=1$ if $k=s_j$ and $\a_k=0$ otherwise. Then $F$ is a transcendental entire function such that $F(\mathfrak{L})\subseteq \mathbb{L}$. In particular, there exist uncountable many transcendental entire functions
taking the set of the ultra-strong Liouville numbers into the set of Liouville numbers.
\end{theorem}

Let us describe in a few words the main ideas for proving Theorem \ref{teo1}. First, our desired function has the form $F(z)=\sum_{n\geq 1}z^{t_n}/10^{t_n!}$, where $(t_n)_n$ is an integer sequence with a very fast growth. We then approximate $F(\xi)$, where $\xi$ is a Liouville number, by a convenient truncation $F_m(p_n/q_n)$ for sufficiently large $m$ and $n$. After that, we take the advantage of the fact that our series has much more zero coefficients than a strongly lacunary series. This, together with the fact that well-approximations come from the continued fraction allows us to arrive at our desired estimate. The proof splits in two cases depending on the growth of the denominator of the approximants of $\xi$.

\section{The proof of Theorem \ref{teo1}}

Let $(s_n)_{n\geq 1}$ and $F(z)$ defined as in the statement of Theorem \ref{teo1}. Clearly, $F$ is a transcendental entire function and now, we shall prove that $F(\mathfrak{L})\subseteq \mathbb{L}$.

Let $\xi$ be a ultra-strong Liouville number and let $(p_n/q_n)_{n\geq 1}$ be the sequence of the convergents of its continued fraction. This means that $0<|\xi-p_n/q_n|<1/q_n^n$, for all $n\geq 1$. Set $\phi_n=\phi_n(\xi)$ as the smallest positive integer $k$ such that $q_n\leq 10^{k!}$. We have two cases to consider:\\

\noindent
{\bf Case 1.} When $\phi_n\leq n^k$ for some $k\geq 1$ and all $n\geq 1$.

In this case, $q_n\leq 10^{\phi_n!}\leq 10^{n^k!}$. Now, consider the truncations
\[
F_n(z):=\sum_{k=0}^{n}\frac{\alpha_k}{10^{k!}}z^k,
\]
and the convergents
\[
\gamma_n:=F_{n}\left(\frac{p_{2n^2}}{q_{2n^2}}\right).
\]
Note that $\den(\gamma_n)=10^{n!}(q_{2n^2})^n$ (where $\den(z)$ denotes the denominator of a rational number $z$). We shall prove that $F(\xi)$ is well-approximated for the rational numbers $\gamma_n$ in a convenient way which ensures that it is a Liouville number. Since $|F(\xi)-\gamma_n| \leq |F(\xi)-F_n(\xi)| + |F_n(\xi)-\gamma_n|$, we need estimate each part in the right-hand side.
For that, we have
\[
F_{n}(\xi)-\gamma_n=\displaystyle\sum_{k=1}^{n}\frac{\a_k}{10^{k!}}\left(\xi^k-\left(\frac{p_{2n^2}}{q_{2n^2}}\right)^k\right).
\]
and it holds that
\[
\left|\xi^k-\left(\frac{p_{2n^2}}{q_{2n^2}}\right)^k\right|\leq \left|\xi-\frac{p_{2n^2}}{q_{2n^2}}\right| \displaystyle\sum_{t=0}^{k-1}|\xi|^{k-j-1}\left|\frac{p_{2n^2}}{q_{2n^2}}\right|^{j}\leq \left|\xi-\frac{p_{2n^2}}{q_{2n^2}}\right|\cdot k(1+|\xi|)^{n-1},
\]
since $\max\{|\xi|,|p_{n}/q_{n}|\}<1+|\xi|$ (for all sufficiently large $n$). Then 
\[
|F_{n}(\xi)-\gamma_n|<\left|\xi-\frac{p_{2n^2}}{q_{2n^2}}\right|\cdot (1+|\xi|)^{n-1}<\frac{(1+|\xi|)^{n-1}}{q_{2n^2}^{2n^2}},
\]
where we used that $\sum_{k\geq 1}k/10^{k!}=0.1200030\ldots$. Since $q_{m}>e^{(m-3)!}$ (which can be proved by using that $p_m/q_m$ is a convergent of the continued fraction of $\xi$), then
\[
q_{2n^2}^{n^2}>e^{n^2(2n^2-3)!}>(1+|\xi|)^{n-1}\cdot 10^{nn!},
\]
for all sufficiently large $n$. Therefore, a straightforward calculation gives
\begin{equation}\label{l3}
|F_{n}(\xi)-\gamma_n|<\frac{1}{(\den(\gamma_n))^n},
\end{equation}
for all sufficiently large $n$.

Now, for estimating $|F(\xi)-F_n(\xi)| $, we shall consider the truncation in $n=s_{j-1}$ satisfying $10^{s_{j-1}}>|\xi|$. Thus, we have
\begin{equation}\label{l4}
|F(\xi)-F_{s_{j-1}}(\xi)| \leq \frac{2}{10^{s_j!-s_{j-1}s_j}}<\frac{1}{(\den(\gamma_{s_{j-1}}))^{s_{j-1}}},
\end{equation}
where we used that $(10^{s_{j-1}!}q_{2s_{j-1}^2}^{s_{j-1}})^{s_{j-1}}\leq 10^{s_j!-s_{j-1}s_j}$ since $s_{j-1}^{4(k+2)}<s_j$ for all sufficiently large $j$ and $q_{2s_{j-1}^2}<10^{(2s_{j-1}^2)^k!}$. By combining (\ref{l3}) and (\ref{l4}), we obtain
\[
|F(\xi)-\gamma_n|<\frac{2}{(\den(\gamma_n))^n},
\]
for all sufficiently large $n$. Thus, in order to prove that $F(\xi)$ is a Liouville number, it suffices to prove that $|F(\xi)-\gamma_n|>0$ for infinitely many integers $n$. Suppose the contrary, then $\gamma_n=p/q$ for all sufficiently large integers $n$. By multiplying this equality by $10^{n!}q_{2n^2}^nq$, we get that $q_{2n^2}$ divides $q$ for infinitely many integers $n$ which is an absurd. Thus $F(\xi)$ is a Liouville number as desired.\\

\noindent
{\bf Case 2.} When $\phi_n$ is not bounded for $n^k$ for all $k\geq 1$.

In this case,  we have the existence of infinitely many pairs $(n_j,k_j)\in \Z_{\geq 1}^2$ such that
\[
\phi_{n_j}\leq n_j^{k_j}\ \mbox{and}\ \phi_{n_j+1}>(n_j+1)^{k_j}.
\]

Now, define $t_j$ as the smallest integer such that $s_{t_j}>\phi_{n_j}$ and define our approximants as
\[
\gamma_{j}:=F_{s_{t_j-1}}\left(\frac{p_{n_j}}{q_{n_j}}\right).
\]
Note that $\den(\gamma_j)=10^{s_{t_j-1}!}q_{n_j}^{s_{t_j-1}}$.

As before, we want to obtain an estimate for  $|F(\xi)-\gamma_j|\leq |F_{s_{t_j-1}}(\xi)-\gamma_j|+|F(\xi)-F_{s_{t_j-1}}(\xi)|$. First, we shall estimate $|F(\xi)-F_{s_{t_j-1}}(\xi)|$. For that, note that for all sufficiently large $j$, we have $10^{s_{t_j-1}}>|\xi|$ and then
\begin{eqnarray*}
  |F(\xi)-F_{s_{t_j-1}}(\xi)|&=&\left|\sum_{k\geq s_{t_j}}\frac{\alpha_k}{10^{k!}}\xi^k\right| \\
           &\leq& \frac{2}{10^{s_{t_j}!-s_{t_j}(s_{t_j-1})}} \\
   &\leq & \frac{1}{10^{s_{t_j-1}!s_{t_j-1}}} \frac{1}{10^{\phi_{n_j}!s_{t_j-1}^2}}\\
&\leq&\frac{1}{10^{s_{t_j-1}!s_{t_j-1}}}\frac{1}{q_{n_j}^{s_{t_j-1}^2}}=\frac{1}{\den(\gamma_j)^{s_{t_j-1}}},
\end{eqnarray*}
where we used that $10^{\phi_n!}\geq q_n$ and $s_{t_j}\geq \min\{s_{t_j-1}^3,\phi_{n_j}+1,5\}$, for all sufficiently large $j$. Then, we have
\begin{equation}\label{l1}
 |F(\xi)-F_{s_{t_j-1}}(\xi)|<\frac{1}{\den(\gamma_j)^{s_{t_j-1}}}.
\end{equation}

Now, we shall estimate $|F_{s_{t_j-1}}(\xi)-\gamma_j|$. For that, we have
\[
F_{s_{t_j-1}}(\xi)-\gamma_j=\displaystyle\sum_{k=1}^{s_{t_j-1}}\frac{\a_k}{10^{k!}}\left(\xi^k-\left(\frac{p_{n_j}}{q_{n_j}}\right)^k\right).
\]
As in the previous case, we get
\[
\left|\xi^k-\left(\frac{p_{n_j}}{q_{n_j}}\right)^k\right|\leq \left|\xi-\frac{p_{n_j}}{q_{n_j}}\right|\cdot k10^{ks_{t_j-1}},
\]
since $\max\{|\xi|,|p_{n_j}/q_{n_j}|\}<1+|\xi|\leq 10^{s_{t_j-1}}$. Then 
\[
|F_{s_{t_j-1}}(\xi)-\gamma_j|<\left|\xi-\frac{p_{n_j}}{q_{n_j}}\right|\cdot 10^{s_{t_j-1}^2}.
\]
Now, we use the well-known fact that
\[
\left|\xi-\frac{p_{n_j}}{q_{n_j}}\right|<\frac{1}{q_{n_j}q_{n_j+1}}.
\]
Also, by definition, since $\phi_{n_j+1}>(n_j+1)^{k_j}$, then $q_{n_j}q_{n_j+1}\geq q_{n_j+1}>10^{(n_j+1)^k!}$. Therefore, for all $k_j\geq 5$, we have
\[
(n_j+1)^{k_j}!\geq (n_j^{k_j}+k_jn_j^{k_j-1})!\geq n_j^{k_j}!k_jn_j^{k_j(k_j-1)/2}\geq 3n_j^{k_j}!n_j^{2k_j}.
\]
Therefore
\begin{equation}\label{l1.5}
|F_{s_{t_j-1}}(\xi)-\gamma_j|<\frac{1}{10^{3n_j^{k_j}!n_j^{2k_j}-s_{t_j-1}^2}}.
\end{equation}
Note that $\den(\gamma_j)^{s_{t_j-1}}\leq 10^{s_{t_j-1}s_{t_j-1}!+\phi_{n_j}!s_{t_j-1}^2}$, since $10^{\phi_{n_j}!}\geq q_{n_j}$. However, we have the inequality
\[
s_{t_j-1}s_{t_j-1}!+\phi_{n_j}!s_{t_j-1}^2\leq 3n_j^{k_j}!n_j^{2k_j}-s_{t_j-1}^2,
\]
since $s_{t_{j}-1}\leq \phi_{n_j}\leq n_j^{k_j}$. The above inequality combined with (\ref{l1.5}) gives
\begin{equation}\label{l2}
|F_{s_{t_j-1}}(\xi)-\gamma_j|<\frac{1}{\den(\gamma_j)^{s_{t_j-1}}}.
\end{equation}
By combining (\ref{l1}) and (\ref{l2}) we obtain
\[
|F(\xi)-\gamma_j|<\frac{2}{\den(\gamma_j)^{s_{t_j-1}}},
\]
for all sufficiently large $j$. Since $|F(\xi)-\gamma_j|>0$ for all sufficiently large $j$ (by a same argument as before), then $F(\xi)$ is a Liouville number as desired. In conclusion, $F(\mathfrak{L})\subseteq \mathbb{L}$. 

The proof of the existence of the uncountable many functions $F(z)$ with this property is because there is a binary tree of different possibilities for $(s_n)_{n\geq 1}$ (and each choice define a different function $F(z)$). For example, take $s_n=a_n^{n!}$, where $a_n\in \{2,3\}$.  
\qed

%Acknowledgements

%\section*{Acknowledgement}
%The author was supported in part by CNPq and FEMAT.

% The Appendices part is started with the command \appendix;
% appendix sections are then done as normal sections
% \appendix

% \section{}
% \label{}

%% BIBLIOGRAPHY

\end{document}